\documentclass[11pt,a4paper]{amsart}
\usepackage{amssymb,amsmath}
\usepackage{mathpazo}

\headsep=1cm \numberwithin{equation}{section}
\newtheorem{theorem}{Theorem}[section]
\newtheorem{lemma}[theorem]{Lemma}
\newtheorem{proposition}[theorem]{Proposition}
\newtheorem{corollary}[theorem]{Corollary}

\theoremstyle{definition}
\newtheorem{definition}[theorem]{Definition}
\theoremstyle{remark}
\newtheorem{remark}[theorem]{Remark}
\newtheorem{example}[theorem]{Example}

\newcommand{\Ass}{\operatorname{Ass}}

\newcommand{\grade}{\operatorname{grade}}
\newcommand{\Assh}{\operatorname{Assh}}
\newcommand{\Spec}{\operatorname{Spec}}

\newcommand{\Ht}{\operatorname{ht}}
\newcommand{\pd}{\operatorname{pd}}
\newcommand{\V}{\operatorname{V}}
\newcommand{\Var}{\operatorname{Var}}

\newcommand{\Ext}{\operatorname{Ext}}
\newcommand{\Supp}{\operatorname{Supp}}

\newcommand{\Hom}{\operatorname{Hom}}

\newcommand{\Ann}{\operatorname{Ann}}

\newcommand{\depth}{\operatorname{depth}}

\newcommand{\fm}{\frak{m}}
\newcommand{\fp}{\frak{p}}
\newcommand{\fq}{\frak{q}}
\newcommand{\fa}{\frak{a}}
\newcommand{\fb}{\frak{b}}

\newcommand{\fn}{\frak{n}}

\begin{document}

\author[Asgharzadeh  and Tousi ]{Mohsen Asgharzadeh
and Massoud Tousi}

\title[Cohen-Macaulayness ... ] {Cohen-Macaulayness with
respect to Serre classes }

\address{M. Asgharzadeh,  Department of Mathematics Shahid Beheshti
University Tehran, Iran, and Institute for Studies in Theoretical
Physics and Mathematics, P.O. Box 19395-5746, Tehran, Iran.}
\email{asgharzadeh@ipm.ir}
\address{M. Tousi, Department of Mathematics Shahid Beheshti
University Tehran, Iran,  and Institute for
Studies in Theoretical Physics and Mathematics, P.O. Box 19395-5746,
Tehran, Iran.} \email{mtousi@ipm.ir} 


\begin{abstract}
Let $R$ be a commutative Noetherian  ring. The notion of regular
sequences with respect to a Serre class of $R$-modules is introduced
and some of their essential properties are given. Then in the local
case, we explore a theory of Cohen-Macaulayness with respect to
Serre classes.
\end{abstract}

\maketitle

\section{Introduction}

The concept of almost vanishing has been studied by a number of
authors, see \cite{GR} and \cite{RSS}. They introduced along the way
the notion of almost zero modules in two different ways over not
necessarily Noetherian rings. We do not require to give the
definition of the almost zero modules, but we list here two basic
properties of them:
\begin{enumerate}
\item[(i)] for any exact
sequence of $A$-modules
$$0\longrightarrow M'\longrightarrow M\longrightarrow M''\longrightarrow
0,$$the module $M$ is almost zero if and only if each of $M'$ and
$M''$ is almost zero, and

\item[(ii)] if $\{M_i\}$ is a directed system consisting of almost zero modules, then
its direct limit $\underset{i} {\varinjlim} M_i$ is almost zero.
\end{enumerate}
A subclass of the class of all modules is called a Serre class, if
it is closed under taking submodules, quotients and extensions. A
serre class which is closed under taking direct limit of any direct
system of its objects is called torsion theory, see \cite{St}.  So,
the class of almost zero modules is a torsion theory.

In view of \cite[Definition 1.2]{RSS}, a sequence
$\underline{x}:=x_{1}, \dots, x_{r}$ of elements of a certain ring
$A$ is called almost regular sequence, if the $A$-module
$((x_{1},\dots, x_{i-1}) :_{A} x_{i})/(x_{1},\dots, x_{i-1})A$ is
almost zero, for each $i = 1, \dots , r$.  If every system of
parameters for A is an almost regular sequence, $A$ is said to be
almost Cohen-Macaulay. We refer the reader to~\cite[Proposition
1.3]{RSS} for a connection between this definition and the monomial
conjecture. These observations motivates us to introduce a new
generalization of the notions of regular sequences and
Cohen-Macaulay modules by using Serre classes. In fact, we do this
by replacing  the class of almost zero modules with an arbitrary
Serre class of modules. It is worth pointing out that some of the
existing generalizations of Cohen-Macaulayness can be viewed as the
special cases of our definition. More precisely, consider the
following two examples.

Let $(R,\fm)$ be a Noetherian local ring and $M$ a finitely
generated $R$-module. Consider the torsion theory
$\mathcal{T}_0:=\{N\in R\textbf{-Mod}|\Supp_R N \subseteq \Supp_R
(R/\fm )\}$, where $R\textbf{-Mod}$ is the category of $R$-modules
and $R$-homomorphisms. In the literature, the concept of
$M$-sequence with respect to $\mathcal{T}_0$ is called $M$-filter
regular sequence, see \cite{CST} and \cite{SV}. Cuong et al.
introduced the notion of f-modules in \cite{CST}. They studied
modules called f-modules satisfying the condition that every system
of parameters is a filter regular sequence, see \cite{CST}. Now,
consider the torsion theory $\mathcal{T}_1:=\{N\in
R\textbf{-Mod}|\dim N\leq 1 \}$. The concept of $M$-sequence with
respect to $\mathcal{T}_1$ is called generalized filter regular
sequence on $M$. This notion, as a generalization of the notion of
filter sequences, was first appeared in \cite{N}. Following
\cite{N}, the concept of generalized f-modules introduced in
\cite{NM}.  An $R$-module $M$ is called a generalized f-module, if
every system of parameters for $M$ is a generalized filter regular
sequence. Thus our theory will include the notions of f-modules and
generalized f-modules.

Throughout this paper, $R$ is a  commutative Noetherian  ring, $\fa$
an ideal of $R$ and $M$ an $R$-module. Always, "$\mathcal{S}$"
stands for a Serre class.

In Section 2, we introduce the notion of weak $M$-sequences with
respect to a Serre class $\mathcal{S}$. We define the
$\mathcal{S}-C.\grade_{R}(\fa,M)$ as the supremum length of weak
$M$-sequences with respect to $\mathcal{S}$ in $\fa$. After
summarizing some results, we characterize
$\mathcal{S}-C.\grade_{R}(\fa,M)$ via some homological tools such as
Ext-modules, Koszul complexes and local cohomology modules. This
provides a common language for expressing some results concerning
several type of sequences that have been appeared in different
papers.

Let $(R,\fm)$ be a Noetherian local ring and $M$ a finitely
generated $R$-module. We say that $M$ is
$\mathcal{S}$-Cohen-Macaulay, if every system of parameters for $M$
is a weak $M$-sequence with respect to $\mathcal{S}$. In Section 3,
we exhibit some of the basic properties of
$\mathcal{S}$-Cohen-Macaulay modules. Some connections between the
notions of $\mathcal{S}$-Cohen-Macaulayness and Cohen-Macaulayness
are given. We show that most of the remarkable properties of
Cohen-Macaulay modules remain valid for $\mathcal{S}$-Cohen-Macaulay
modules. Especially, they behave well with respect to flat
extensions, annihilators of local cohomology modules, non
Cohen-Macaulay locus and quotient by weak sequences with respect to
the Serre class $\mathcal{S}$.

\section{Grade of ideals with respect to Serre classes}

Let $\mathcal{S}$ be a subcategory of the category of $R$-modules
and $R$-homomorphisms. Then $\mathcal{S}$ is said to be a
\textit{Serre class} (or Serre subcategory), if for any exact
sequence of $R$-modules
$$0\longrightarrow L\longrightarrow M\longrightarrow N\longrightarrow
0,$$ the $R$-module $M$ belongs to $\mathcal{S}$ if and only if each
of $L$ and $N$ belongs to $\mathcal{S}$. In the case $\mathcal{S}$
is subclass of finitely generated $R$-modules, the following may be
useful.

$\mathcal{S}$ is Serre class if and only if the following conditions
are satisfied by $\mathcal{S}$:
\begin{enumerate}
\item[(i)]
the kernel (resp.  the cokernel) of every morphism of objects of
$\mathcal{S}$ is also in $\mathcal{S}$, and
\item[(ii)] for any exact
sequence $0\longrightarrow L\longrightarrow M\longrightarrow
N\longrightarrow 0$, if both $L$ and $N$ are objects of
$\mathcal{S}$, then so is $M$.
\end{enumerate}
For the proof see, \cite[Corollary 3.2]{T}.

The key to the work in this paper is given by the following easy
lemma.

\begin{lemma}\label{key}
Let $M$ be a finitely generated $R$-module and $K$ an $R$-module.
Then the following hold:
\begin{enumerate}
\item[(i)] $M\in\mathcal{S}$
if and only if  $R/ \fp\in\mathcal{S}$, for all $\fp \in \Supp M$.
In particular, for any two finitely generated $R$-modules $N,L$ with
$\Supp_R N=\Supp_R L$, we have $N\in\mathcal{S}$ if and only if
$L\in\mathcal{S}$.
\item[(ii)]If $S$ is closed under taking direct sums,
then $K\in\mathcal{S}$ if and only if  $R/ \fp\in\mathcal{S}$, for
all $\fp \in \Supp K$. In particular, for any two $R$-modules $N,L$
with $\Supp_R N=\Supp_R L$, it follows that $N\in\mathcal{S}$ if and
only if $L\in\mathcal{S}$.
\end{enumerate}
\end{lemma}

{\bf Proof.} (i) Without loss of generality we can assume that
$M\neq0$. First, assume that $M\in\mathcal{S}$. Let $\fp\in
\Supp_{R}M$. So there exists $m\in M$ such that $(0:_{R}m)\subseteq
\fp$. Since $R/(0:_{R}m)\cong Rm \subseteq M$, it turns out that
$R/(0:_{R}m)\in\mathcal{S}$. From the natural epimorphism
$R/(0:_{R}m)\longrightarrow R/\fp$, we get $R/\fp\in\mathcal{S}$.

Now we prove the converse. There is a chain $$0=M_{0}\subseteq
M_{1}\subseteq\dots\subseteq M_{\ell}=M$$ of submodules of $M$ such
that for each $j$, $M_{j}/M_{j-1}\cong R/\fp$ for some $\fp\in \Supp
M$. By using short exact sequences the situation can be reduced to
the trivial case $\ell = 1$.

(ii) Since $\mathcal{S}$ is closed under taking direct limits, we
can assume that $K$ is a finitely generated $R$-module. The
remainder of the proof is similar to (i). $\Box$

\begin{definition}Let $M$ be an $R$-module. A sequence
$\underline{x}:=x_{1}, \dots, x_{r}$ of elements of $R$ is called a
weak $M$-sequence with respect to $\mathcal{S}$ if for each $i = 1,
\dots , r$ the $R$-module $((x_{1},\dots, x_{i-1}) :_{M}
x_{i})/(x_{1},\dots, x_{i-1})M$ belongs to $\mathcal{S}$. If in
addition $M/\underline{x} M\notin\mathcal{S}$, we say that
$\underline{x}$ is an $M$-sequence with  respect to $\mathcal{S}$.
\end{definition}

\textbf{Notation.} For an $R$-module $L$, we denote $\{\fp \in
\Supp_R L |R/ \fp \notin\mathcal{S}\}$ by $\mathcal{S}-\Supp_{R}L$
and $\{\fp \in \Ass_R L |R/\fp \notin\mathcal{S}\}$ by
$\mathcal{S}-\Ass_{R}L$.\\

In order to exploit Lemma \ref{key} and Definition 2.2, we give the
relation between regular sequences with respect to $\mathcal{S}$ and
ordinary  regular sequences.

\begin{lemma} Let $M$ be a finitely generated $R$-module and
$\underline{x}:=x_{1}, \dots, x_{r}$ a sequence of elements of $R$.
Then the following conditions are equivalent:
\begin{enumerate}
\item[(i)] $x_{i}\notin \bigcup_{\fp\in
\mathcal{S}-\Ass_{R}M/(x_{1},\cdots,x_{i-1})M}\fp$ for all $i=1,
\cdots, r$.
\item[(ii)] The sequence  $x_{1}, \dots, x_{r}$ is a weak $M$-sequence with
respect to $\mathcal{S}$.
\item[(iii)] For any $\fp \in \mathcal{S}-\Supp_{R}(M)$, the elements
$x_{1}/1, \dots, x_{r}/1 $  of the local ring $R_{\fp}$ form a weak
$M_{\fp}$-sequence.
\item[(iv)] The sequence $x^{n_1}_1
,\cdots, x^{n_r}_r$ is a weak $M$-sequence with respect to
$\mathcal{S}$ for all positive integers $n_1,\cdots,n_r$.
\end{enumerate}
\end{lemma}

{\bf Proof.} $(i)\Rightarrow (ii)$ Let $1\leq i\leq r$. In view of
Lemma \ref{key}, it is enough to show that $\{R/\fp:\fp \in
\Supp_R((x_{1},\dots, x_{i-1}) M:_{M} x_{i})/(x_{1},\dots,
x_{i-1})M\}\subseteq \mathcal{S}$. To establish this, suppose on the
contrary that, there is $\fp \in\Supp_R((x_{1},\dots, x_{i-1})
M:_{M} x_{i})/(x_{1},\dots, x_{i-1})M)$ such that $R/ \fp
\notin\mathcal{S}$. So there exists $\fq\in \Ass_R ((x_{1},\dots,
x_{i-1}) :_{M} x_{i})/(x_{1},\dots, x_{i-1})M$, which is contained
in $\fp$. The natural epimorphism $R/\fq\longrightarrow R/\fp$ and
the condition $R/\fp\notin\mathcal{S}$ imply that
$R/\fq\notin\mathcal{S}$. Since  $\fq\in \mathcal{S}-\Ass_{R}(M/
(x_{1},\dots, x_{i-1})M)$ and $x_i \in \fq$, we get a contradiction.

$(ii)\Rightarrow (iii)$ This implication is an immediate consequence
of Lemma \ref{key}.

$(iii)\Rightarrow (i)$  Assume that
$\fp\in\mathcal{S}-\Ass_{R}(M/(x_{1}, \dots, x_{i-1})M)$, for some
$i = 1, \dots , r$. Then $\fp\in \mathcal{S}-\Supp_{R}(M)$ and $\fp
R_{\fp}\in\Ass_{R_{\fp}}(M_{\fp}/(x_{1}/1, \dots,
x_{i-1}/1)M_{\fp})$. Hence by our assumptions we have $x_{i}/1\notin
\fp R_{\fp}$. Consequently $x_{i}\notin \fp$.

$(iii)\Leftrightarrow (iv)$ is clear.
$\Box$\\

Now, we establish a preliminary lemma.

\begin{lemma} \label{ext}
Let $\underline{x}:=x_{1}, \dots, x_{r}$ be a weak $M$-sequence with
respect  to $\mathcal{S}$ in $\fa$. Then the following hold:
\begin{enumerate}
\item[(i)] $\Ext^{i}_R(R/\fa,M)\in \mathcal{S}$ for all $0\leq i\leq
r-1$.
\item[(ii)]$\Ext^{r}_R(R/\fa,M)\notin \mathcal{S}$ if and
only if $\Hom_R(R/\fa,M/\underline{x}M)\notin \mathcal{S}$.
\end{enumerate}
\end{lemma}

{\bf Proof.} (i) Let $0\leq i\leq r-1$ and $\fp \in
\Supp_R(\Ext^{i}_R(R/\fa,M))$.  Assume that $R/ \fp
\notin\mathcal{S}$. We have $\fp \in \Supp_R(M)$. By the implication
$ii)\Rightarrow iii)$ of Lemma 2.3, the elements $x_{1}/1, \dots,
x_{r}/1 $  of the local ring $R_{\fp}$ form a weak
$M_{\fp}$-sequence. Therefore
$\Ext^{i}_{R}(R/\fa,M)_{\fp}\cong\Ext^{i}_{R_{\fp}}(R_{\fp}/\fa_{\fp},
M_{\fp})=0$, which is a contradiction. So $R/ \fp \in\mathcal{S}$.
Now Lemma \ref{key} implies that $\Ext^{i}_R(R/\fa,M)\in
\mathcal{S}$.

(ii) Let $\fp \in \mathcal{S}-\Supp_{R}(M)$. By Lemma 2.3, the
elements $x_{1}/1, \dots, x_{r}/1 $ of the local ring $R_{\fp}$ form
a weak $M_{\fp}$-sequence. Therefore
$\Ext^{r}_{R_{\fp}}(\frac{R_{\fp}}{\fa
R_{\fp}},M_{\fp})\cong\Hom_{R_{\fp}}(R_{\fp}/\fa
R_{\fp},M_{\fp}/\underline{x}M_{\fp})$. This shows that
$\mathcal{S}-\Supp_{R}(\Ext^{r}_R(R/\fa,M))
=\mathcal{S}-\Supp_{R}(\Hom_R(R/\fa,M/\underline{x}M))$. Hence the
desired result follows from Lemma \ref{key}. $\Box$

Let $\fa$ be an ideal of the ring $R$. Suppose  $\underline{y}:=
y_{1}, \cdots, y_{r}$ is a system of generators of $\fa$. We denote
the Koszul complex of $\underline{y}$ by
$K_{\bullet}(\underline{y})$. For an R-module $M$, the  Koszul
complex with coefficient in $M$, is defined by
$K^{\bullet}(\underline{y}, M) :=\Hom( K_{\bullet}(\underline{y}),
M)$.

We need the following lemma in Definition 2.6 below. The symbol
$\mathbb{N}_{0}$ will denote the set of non-negative integers.

\begin{lemma}
Let $\fa$ be an ideal of the ring $R$ and $M$ an $R$-module. Suppose
that $\underline{y}:=y_{1}, \cdots, y_{r}$ is a system of generators
of $\fa$. The notion $\inf\{i\in
\mathbb{N}_{0}|H^{i}(K^{\bullet}(\underline{y},M))\notin
\mathcal{S}\}$  does not depend on the choice of the generating sets
of $\fa$.
\end{lemma}

{\bf Proof.} Let $\underline{x}:= x_{1}, \cdots, x_{s}$ be another
generating set for $\fa$. Set $\underline{x}':= x_{1}, \cdots,
x_{s},y_{1}, \cdots, y_{r}$. In view of \cite[Proposition
1.6.21]{BH}, $H_{\bullet}(\underline{x}',M)\cong
H_{\bullet}(\underline{x},M)\otimes_R \bigwedge R^{r}$. Therefore
$H_{i}(\underline{x}',M)\in\mathcal{S}$ if and only if
$H_{i}(\underline{x},M)\in\mathcal{S}$, since $\bigwedge R^{r}$ is a
finitely generated free $R$-module. Set $n=r+s$. Thus the symmetry
of Koszul cohomology and Koszul homology implies that
$H^{n-i}(K^{\bullet}(\underline{x}',M))\in\mathcal{S}$ if and only
if $H^{n-i}(K^{\bullet}(\underline{x},M))\in\mathcal{S}$, see
\cite[Proposition 1.6.10(d)]{BH}. Set $\underline{y}':=y_{1},
\cdots, y_{r}, x_{1}, \cdots, x_{s}$. By the same reason we have
$H^{n-i}(K^{\bullet}(\underline{y}',M))\in\mathcal{S}$ if and only
$H^{n-i}(K^{\bullet}(\underline{y},M))\in\mathcal{S}$. The claim
becomes clear from  the fact that Koszul complex
$K^{\bullet}(\underline{y}',M)$ is invariant (up to isomorphism)
under permutation of  $\underline{y}'$. $\Box$

\begin{definition}Let  $M$ be an $R$-module and $\fa$ an ideal of
$R$.  Let $\underline{x}:= x_{1}, \cdots, x_{s}$ be a generating set
for $\fa$. The notions of local cohomology grade, Ext grade, Koszul
grade and classical grade of $\fa$ on $M$ with respect to
$\mathcal{S}$, are defined, respectively as follows:
\begin{enumerate}
\item[(i)]$\mathcal{S}-H.\grade_{R}(\fa,M):=\inf\{i\in
\mathbb{N}_{0}|H_{\fa}^{i}(M)\notin \mathcal{S}\}$,
\item[(ii)]$\mathcal{S}-E.\grade_{R}(\fa,M):=\inf\{i\in
\mathbb{N}_{0}|\Ext^{i}_R(R/\fa,N)\notin \mathcal{S}\}$,
\item[(iii)]
$\mathcal{S}-K.\grade_{R}(\fa,M):=\inf\{i\in
\mathbb{N}_{0}|H^{i}(K^{\bullet}(\underline{x},M))\notin
\mathcal{S}\}$,
\item[(iv)]
$\mathcal{S}-C.\grade_{R}(\fa,M):=\sup \{\ell\in \mathbb{N}_{0}|
y_{1}, \cdots, y_{\ell}  \emph{ is a weak M-sequence  in } \fa
\emph{ with}\\ \emph{ respect to } \mathcal{S} \}$.
\end{enumerate} Here $\inf$ and $\sup$
are formed in $\mathbb{Z} \cup \{\pm\infty\}$ with the convention
that $\inf \emptyset=+ \infty$ and $\sup \emptyset=-\infty$.
\end{definition}

In the case $\mathcal{S}=\{0\}$, for simplicity we use the notions
$K.\grade_R(\fa,M)$, $E.\grade_R(\fa,M)$ and $H.\grade_R(\fa,M)$,
instead of $\mathcal{S}-K.\grade_{R}(\fa,M)$,
$\mathcal{S}-E.\grade_{R}(\fa,M)$ and
$\mathcal{S}-H.\grade_{R}(\fa,M)$.

\begin{proposition}\label{local}Let $M$ be a finitely generated
$R$-module. Then  the following hold:
\begin{enumerate}
\item[(i)]$\mathcal{S}-E.\grade_{R}(\fa,M)=\inf\{E.\grade
_{R_{\fp}}(\fa R_{\fp},M_{\fp})
|\fp\in \mathcal{S}-\Supp_{R}(M)\}$.
\item[(ii)]
$\mathcal{S}-K.\grade_{R}(\fa,M)=\inf\{K.\grade_{R_{\fp}}(\fa
R_{\fp},M_{\fp})|\fp\in\mathcal{S}-\Supp_{R}(M)\}$.
\item[(iii)] If $\mathcal{S}$
is closed under taking direct sums,
then$$\mathcal{S}-H.\grade_{R}(\fa,M)=\inf\{H.\grade_{R_{\fp}}(\fa
R_{\fp},M_{\fp})|\fp\in\mathcal{S}-\Supp_{R}(M)\}.$$
\end{enumerate}
\end{proposition}

{\bf Proof.} (i) First assume that
$s:=\mathcal{S}-E.\grade_{R}(\fa,M)<\infty$.  So $\Ext^{s}_R
(R/\fa,M)\notin\mathcal{S}$ and $\Ext^{i}_R
(\frac{R}{\fa},M)\in\mathcal{S}$ for all $0\leq i<s$. By using Lemma
\ref{key}(i), $\Ext^{i}_{R_{\fp}}(\frac{R_{\fp}}{\fa R_{\fp}}
,M_{\fp})=0$ for all $\fp\in \mathcal{S}-\Supp_{R}(M)$ and all
$0\leq i<s$. Therefore $s\leq E.\grade_{R_{\fp}}(\fa
R_{\fp},M_{\fp})$ for all $\fp\in \mathcal{S}-\Supp_{R}(M)$. Again
by Lemma \ref{key}(i), $R/\fq\notin\mathcal{S}$ for some
$\fq\in\Supp_R(\Ext^{s}_R(\frac{R}{\fa},M))$. The fact
$\Ext^{s}_{R_{\fq}}(\frac{R_{\fq}}{\fa R_{\fq}},M_{\fq})\neq0$
implies that $E.\grade_{R_{\fq}}(\fa R_{\fq},M_{\fq})\leq s$.

Now, assume that $\mathcal{S}-E.\grade_{R}(\fa,M)=\infty$.  In the
case $\mathcal{S}-\Supp_{R}(M)=\emptyset$, we have nothing to prove.
Hence we can assume that $\mathcal{S}-\Supp_{R}(M)\neq\emptyset$.
Let $\fp\in\mathcal{S}-\Supp_{R}(M)$. Since $\Ext^{i}_R
(R/\fa,M)\in\mathcal{S}$ for all $i\geq 0$, so
$\Ext^{i}_{R_{\fp}}(R_{\fp}/\fa R_{\fp} ,M_{\fp})=0$ for all $i\geq
0$. Consequently $E.\grade_{R_{\fp}}(\fa R_{\fp},M_{\fp})=\infty$.

By the same argument as (i), we can prove (ii) and (iii). Only note
that in the case $(iii)$ we use Lemma \ref{key}(ii) instead of Lemma
\ref{key}(i), since $H^i _{\fa}(M)$ is not necessary finitely
generated $R$-module. $\Box$\\

The following is the main result of this section.

\begin{theorem}\label{comm}
Let $M$ be a finitely generated $R$-module. Then the following hold:
\begin{enumerate}
\item[(i)] Let $\underline{x}:=x_{1},
\dots, x_{r}$ be a maximal weak $M$-sequence with respect  to
$\mathcal{S}$ in $\fa$. Then $r=\mathcal{S}-E.\grade_{R}(\fa,M)$.

\item[(ii)]$\mathcal{S}-C.\grade_{R}(\fa,M)=
\mathcal{S}-E.\grade_{R}(\fa,M)=\mathcal{S}-K.\grade_{R}(\fa,M)$.
\item[(iii)]If $\mathcal{S}$
is closed under taking direct sums, then
$\mathcal{S}-E.\grade_{R}(\fa,M)= \mathcal{S}-H.\grade_{R}(\fa,M)$.

\item[(iv)] If $M/\fa M\notin\mathcal{S}$, then
$\mathcal{S}-K.\grade_R(\fa,M)<\infty$ and all maximal M-sequence
with respect to $\mathcal{S}$ in $\fa$ have a common length.
\end{enumerate}
\end{theorem}

{\bf Proof.} (i) Since $M$ is a finitely generated, it follows from
the maximality of $\underline{x}$  and Lemma 2.3 that $\fa \subseteq
\fp$ for some $\fp\in \mathcal{S}-\Ass_{R}(M/\underline{x}M)$. Hence
$$\fp\in \Ass_R (M/ \underline{x}M)\cap\Supp_R (R/ \fa )= \Ass
_R(\Hom_{R}(R/ \fa,M/ \underline{x}M)),$$ see \cite[Exercise
1.2.27]{BH}. Lemma \ref{key} implies that
$\Hom_{R}(R/{\fa},M/\underline{x}M)\notin\mathcal{S}$. Now, the
conclusion  follows from Lemma \ref{ext}.

(ii) The inequality  $\mathcal{S}-C.\grade_{R}(\fa,M)\leq
\mathcal{S}-E.\grade_{R}(\fa,M)$, becomes clear by Lemma 2.4(i).
Therefore without loss of generality we can assume that
$r:=\mathcal{S}-C.\grade_{R}(\fa,M)<\infty$. So, the other side
inequality follows from (i).

In order to prove the second equality, recall that for all $\fp\in
\mathcal{S}-\Supp_{R}(M)$ we have $E.\grade_{R_{\fp}}(\fa
R_{\fp},M_{\fp})=K.\grade_{R_{\fp}}(\fa R_{\fp},M_{\fp})$ see
\cite[Theorem 6.1.6]{Str}. In view  of Proposition \ref{local}(i)
and (ii) the assertion follows.

(iii) This follows from Proposition \ref{local}(iii) and
\cite[Proposition 5.3.15]{Str}, which state that
$E.\grade_{R_{\fp}}(\fa R_{\fp},M_{\fp})=H.\grade_{R_{\fp}}(\fa
R_{\fp},M_{\fp})$.

(iv) Assume that $\fa$ can be generated by $n$ elements
$\underline{y}:=y_{1}, \cdots, y_{r}$. We have
$H^{n}(K(\underline{y},M))\cong M/ \fa M \notin \mathcal{S}$ and
consequently $\mathcal{S}-K.\grade_R(\fa,M)<\infty$. So, in view of
(i) and (ii) all maximal M-sequences with respect to $\mathcal{S}$
in $\fa$ have the same length. $\Box$

\begin{example}(i)
In  Theorem \ref{comm}(iii) the assumption "$\mathcal{S}$ is closed
under taking direct sums" is really need. To see this, let $(R,\fm)$
be a Cohen-Macaulay local ring of dimension $d>0$ and $\mathcal{S}$
the Serre class of all finitely generated $R$-modules. It is easy to
see that $\mathcal{S}-E.\grade_{R}(\fm,R) =\infty\neq
d=\mathcal{S}-H.\grade_{R}(\fm,R)$.

(ii) In  Theorem \ref{comm}(ii) the finitely generated assumption on
$M$ is necessary. To see this, let $R=K[[x,y]]$ and set
$M:=\bigoplus_{0\neq r\in(X,Y)} R/rR$. By \cite[Page 91]{Str}, we
have $E.\grade_R(\fm,M)=1$ and $C.\grade_R(\fm,M)=0$.
\end{example}

We denote the category of $R$-modules and $R$-homomorphisms, by
$R\textbf{-Mod}$. For an $R$-module $L$, set $\mathcal{T}(L):=\{N\in
R\textbf{-Mod}|\Supp_R N\subseteq \Supp_R L\}$. It is easy to see
that $\mathcal{T}(L)$ is a Serre class, which is closed under taking
direct sums. Such Serre classes are called torsion theories. In the
following proposition we give a characterization of torsion theories
over Noetherian rings.

\begin{proposition}
Let $\mathcal{T}$ be a torsion theory. Then
$\mathcal{T}=\mathcal{T}(L)$, for some $R$-module $L$.
\end{proposition}

{\bf Proof.}  Without loss of generality we  can assume that
$\mathcal{T}\neq0$. If not, we can take $L=0$. Set
$L:=\bigoplus_{\fp\in\sum}  R/\fp $, where $\sum=\{\fp\in \Spec R |
R/\fp\in\mathcal{T}\}$ is a non empty subset of $\Spec R $. Let $N$
be an  $R$-module. It follows from Lemma \ref{key}(ii) that
$N\in\mathcal{T}$ if and only if $\Supp_R N \subseteq \sum$. The
claim follows from the fact that $\Supp_R N\subseteq \sum$ if and
only if $N\in\mathcal{T}(L)$. $\Box$

Let $(R,\fm)$ be a local ring and set $\mathcal{T}(R/\fm):=\{N\in
R\textbf{-Mod}|\Supp_R N \subseteq \Supp_R (R/\fm )\}$. In the
literature, the concept of  $M$-sequence with respect to
$\mathcal{T}(R/\fm)$ is called  $M$-filter regular sequence and  the
maximal length of such sequences in $\fa$ is denoted by
$f-\depth_{R}(\fa,M)$, see \cite{CST} and \cite{Mel2}. The following
corollary can be found in \cite[Proposition 3.2]{LT},  \cite[Theorem
3.9]{LT}, \cite[Theorem 3.10]{LT} and \cite[Theorem 5.5]{Mel1}.

\begin{corollary}
Let $(R,\fm)$ be a local ring and $M$ a finitely generated
$R$-module such that $\Supp_R (M/\fa M) \nsubseteq \{\fm\}$. Then

\[\begin{array}{ll} f-\depth_{R}(\fa,M)&= \inf\{i : H^i _\fa(M)
\textit{ is not Artinian}\}\\&=\inf\{i : H^i (K^\bullet(\fa,M))
\textit{ is not Artinian}\}\\&=\inf\{i : \Ext^i _R(R/\fa,M) \textit{
is not Artinian}\}
.\\
\end{array}\]
\end{corollary}

{\bf Proof.} Note that $\inf\{i : H^i _\fa(M) \textit{ is not
Artinian}\}= \inf\{i : \Supp_R H^i _\fa(M) \nsubseteq \{\fm\}\}$,
see \cite[Lemma 2.4]{Ma}. Also we know that for any finitely
generated $R$-module $N$, $\Supp_R N\subseteq \{\fm\}$ if and only
if $N$ is an Artinian $R$-module. The desired result follows from
Theorem \ref{comm}. $\Box$

By \cite[Definition \ref{key}]{N}, a sequence
$\underline{x}:=x_1,\cdots,x_r$ of elements  of $\fa$ is  a
generalized regular sequence of $M$, if $x_i\notin \fp$ for all
$\fp\in \Ass_R (M/(x_1,\cdots, x_{i-1})M)$ satisfying $\dim (R/\fp)
> 1$. She called, the length of a maximal generalized
regular sequence of $M$ in $\fa$, generalized depth of M in $\fa$
and denoted it by $g-\depth_R(\fa,M)$, see \cite[Definition 4.2]{N}.
Set $X_1:=\{\fp\in \Spec R:\dim (R/\fp)\leq 1\}$ and
$M_1:=\bigoplus_{\fp\in X_1}R/\fp$. Consider the
$\mathcal{T}:=\mathcal{T}(M_1)=\{N\in R\textbf{-Mod}|\Supp_R N
\subseteq \Supp_R  M_1  \}=\{N\in R\textbf{-Mod}|\dim N\leq 1\}$. It
is easy to see that $\underline{x}$ is  a generalized regular
sequence of $M$ if and only if $\underline{x}$ is an $M$-sequence
with respect to $\mathcal{T}$. In the local case, the first and the
second equality in the following corollary is in \cite[Proposition
4.4]{N} and \cite[Proposition 4.5]{N}.

\begin{corollary}
Let $M$ be a finitely generated $R$-module such that $\dim (M/\fa
M)\geq 2$. Then
\[\begin{array}{ll} g-\depth_R(\fa,M)&=
\inf\{i : \dim H^i _\fa(M)
>1\}\\&=\inf\{i : \dim \Ext^i _R (R/\fa,M)>1\}\\&=\inf\{i : \dim H^i (K^\bullet(\fa,M))
>1\}.\\
\end{array}\]
\end{corollary}

Let $j$ be an integer such that $0\leq j <\dim R$ and set
$X_j=\{\fp\in \Spec R:\dim (R/\fp)\leq j\}$. Consider the $R$-module
$M_j:=\bigoplus_{\fp\in X_j}R/\fp$ and set
$$\mathcal{T}_j:=\mathcal{T}(M_j)=\{N\in R\textbf{-Mod}|\Supp_R N
\subseteq \Supp_R  M_j \}=\{N\in R\textbf{-Mod}|\dim N\leq j\}.$$ In
the local case, the first equality in the following corollary was
first appeared in \cite[Page 9]{Q} and \cite[Lemma 2.4]{BN}. The
second equality is in \cite[Lemma 2.3]{CHK}.

\begin{corollary}
Let $M$ be a finitely generated $R$-module such that $\dim (M/\fa
M)\geq j+1$. Then
\[\begin{array}{ll}
\mathcal{T}_j -\depth_R(\fa,M)&=\inf\{i : \dim H^i _\fa(M)
>j\}\\&=\inf\{i : \dim \Ext^i
_R(R/\fa,M)>j\}\\&=\inf\{i : \dim H^i (K(\fa,M))>j\}.\\
\end{array}\]
\end{corollary}

Let $\fb$ be an ideal of $R$. Recall  from \cite[Definition 1.3]{A}
that a sequence $a_1,\cdots ,a_r$  is  a $\fb$-filter regular
$M$-sequence if $x_i\notin \fp$ for all $\fp\in \Ass_{R}
M/(x_1,\cdots, x_{i-1})M\setminus \V(\fb)$. Consider the torsion
theory $\mathcal{T}_\fb:=\mathcal{T}(R/\fb)=\{N\in
R\textbf{-Mod}|\Supp_R N \subseteq \Supp_R (R/\fb)\}$. It is easy to
see that $\underline{x}$ is $\fb$-filter regular $M$-sequence if and
only if $\underline{x}$ is $M$-sequence with respect to
$\mathcal{T}_\fb$. The first equality in the following corollary was
first appeared in \cite[Theorem 1.7]{A}.

\begin{corollary}
Let $M$ be a finitely generated $R$-module such that $\Supp_R (M/\fa
M) \nsubseteq \V(\fb)$. Then
\[\begin{array}{ll}
\mathcal{T}_\fb-\grade(\fa,M)&= \inf\{i :  \Ext^i _R
(R/\fa,M)\notin\mathcal{T}_\fb\}\\&=\inf\{i :  H^i
(K^\bullet(\fa,M))\notin\mathcal{T}_\fb\}\\&=\inf\{i :  H^i _\fa(M)
\notin\mathcal{T}_\fb\}.\\
\end{array}\]
\end{corollary}

\section{Cohen-Macaulayness with respect to Serre classes}
In this section we introduce the concept of
$\mathcal{S}$-Cohen-Macaulay modules. First of all, consider the
following definition.

\begin{definition}Let $M$ be an $R$-module and $\fa$ an ideal of
$R$. The height of $\fa$ on $M$ and the Krull dimension of $M$ with
respect to $\mathcal{S}$ are defined as follows:
\begin{enumerate}
\item[(i)]
$\mathcal{S}-\Ht_{M}(\fa):=\inf\{\Ht_{M}(\fq)|\fq\in
 \mathcal{S}-\Supp_{R}(M)\cap\V(\fa)\}$,

\item[(ii)]$\mathcal{S}-\dim (M)=\sup\{\Ht_{M}(\fq)|\fq\in
 \mathcal{S}-\Supp_{R}(M)\}$.
\end{enumerate}
\end{definition}
In the following lemma we investigate the relation between
$\mathcal{S}-E.\grade_{R}(\fa,M)$ and $\mathcal{S}-\Ht_{M}(\fa)$.

\begin{lemma}
Let $M$  be a finitely generated $R$-module  and $\fa$ an ideal of
$R$. Then $\mathcal{S}-E.\grade_{R}(\fa,M)\leq
\mathcal{S}-\Ht_{M}(\fa)$.
\end{lemma}

{\bf Proof.} The result follows from the following inequality and
equalities.
\[\begin{array}{ll}
\mathcal{S}-E.\grade_{R}(\fa,M)&=\inf\{E.\grade_{R_{\fp}}(\fa
R_{\fp},M_{\fp}) |\fp\in
\mathcal{S}-\Supp_{R}(M)\}\\&=\inf\{E.\grade_{R_{\fp}}(\fa
R_{\fp},M_{\fp}) |\fp\in \mathcal{S}-\Supp_{R}(M/\fa M)\}
\\&\leq\inf\{\Ht_{M_{\fp}}(\fa R_{\fp})
|\fp\in \mathcal{S}-\Supp_{R}(M/\fa M)\}\\&=\inf\{\Ht_{M_{\fp}}(\fp
R_{\fp}) |\fp\in \mathcal{S}-\Supp_{R}(M/\fa
M)\}\\&=\mathcal{S}-\Ht_{M}(\fa),
\\
\end{array}\]
where the first equality follows from the Proposition \ref{local}.
$\Box$

For a subset $X$ of $\Spec R$, we denote  the set of minimal members
of $X$ with respect to inclusion, by $\min (X)$. We say that an
ideal $\fa$ of $R$ is  unmixed on $M$ with respect to $\mathcal{S}$,
if $\mathcal{S}-\Ass_{R}(M/\fa M)=\{\fp\in\min(\Supp_{R}(M/ \fa
M)):R/ \fp\notin\mathcal{S}\}$.

\begin{proposition}\label{localization}
Let $M$ be a finitely generated $R$-module. Then the following are
equivalent:
\begin{enumerate}\item[(i)]$\mathcal{S}-E.\grade _{R}(\fa,M)=
\mathcal{S}-\Ht_{M}(\fa)$ for all ideal $\fa$ of $R$.
\item[(ii)]  $\mathcal{S}-E.\grade_{R}(\fp,M)=\mathcal{S}-\Ht_{M}
(\fp)=\Ht_{M}(\fp)$ for all $\fp\in\mathcal{S}-\Supp_{R}(M)$.
\item[(iii)] For any $\fp\in\mathcal{S}-\Supp_{R}(M)$, $M_{\fp}$ is
Cohen-Macaulay.
\item[(iv)]  Any ideal $\fa$ which generated by
$\Ht_M (\fa)$ elements is unmixed on $M$ with respect to
$\mathcal{S}$.
\end{enumerate}
\end{proposition}

{\bf Proof.} $(i)\Rightarrow (ii)$ Note that
$\mathcal{S}-\Ht_{M}(\fp)=\Ht_{M}(\fp)$ for each
$\fp\in\mathcal{S}-\Supp_{R}(M)$. So this implication is clear.

$(ii)\Rightarrow (iii)$ Let $\fp\in\mathcal{S}-\Supp_{R}(M)$. Then
the claim follows from the following inequalities:

$$
\begin{array}{ll}
\mathcal{S}-E.\grade_{R}(\fp,M)&\leq E-\grade_{R_{\fp}}(\fp
R_{\fp},M_{\fp})\\&\leq\Ht_{M_{\fp}}(\fp
R_{\fp})\\&=\mathcal{S}-\Ht_{M}(\fp)\\&=\mathcal{S}-E.\grade_{R}(\fp,M),
\end{array}
$$

where the first inequality follows from the Proposition \ref{local}.

$(iii)\Rightarrow (iv)$ Let $\fa$  be an ideal, which can be
generated by $\Ht_M (\fa)$ elements, $\{a_1,\cdots,a_n\}$. Let
$\fp\in\mathcal{S}-\Ass_{R}(M/\fa M)$. Then  $M_{\fp}$ is
Cohen-Macaulay and $\Ht_{M_{\fp}}(\fa R_{\fp})=n$. Therefore the
elements $a_{1}/1, \dots, a_{n}/1 $  of the local ring $R_{\fp}$
form an $M_{\fp}$-sequence. So $M_{\fp}/ \fa M_{\fp}$ is
Cohen-Macaulay. Since $\fp R_{\fp}\in\Ass_{R_{\fp}} (M_{\fp}/\fa
M_{\fp})$, we get that $\fp\in\min(\Supp_{R}(M/ \fa M))$.

$(iv)\Rightarrow (iii)$ Let $\fp\in\mathcal{S}-\Supp_{R}(M)$ be such
that $\Ht_M (\fp)=n$. From this one can find the elements
$x_1,\cdots,x_n$ of the ideal $\fp$ such that
$\Ht_M(x_1,\cdots,x_i)=i$ for all $1\leq i\leq n$. By our
assumption, $x_{i}\notin \bigcup_{\fq\in
\mathcal{S}-\Ass_{R}\frac{M}{(x_{1},\cdots,x_{i-1})M}}\fq$, for all
$1\leq i\leq n$. In view of Lemma 2.3, the elements $x_{1}/1, \dots,
x_{n}/1$  of the local ring $R_{\fp}$ become an $M_{\fp}$-sequence.
Therefore $\depth_{R_{\fp}}(M_{\fp})\geq \dim M_{\fp}$, which is we
want.

$(iii)\Rightarrow (i)$ If $\mathcal{S}-\Supp_{R}(M)=\emptyset$, we
have nothing to prove. So we can assume that
$\mathcal{S}-\Supp_{R}(M)\neq\emptyset$. The proof in this case
follows from the proof of Lemma 3.2. $\Box$

\begin{definition} Let $(R,\fm)$ be a local ring and $M$ a
finitely generated $R$-module. Then $M$ is called  an
$\mathcal{S}$-\textit{Cohen-Macaulay} $R$-module if any system of
parameters of M form a weak $M$-sequence with respect to
$\mathcal{S}$. The  ring $R$ is called $\mathcal{S}$-Cohen-Macaulay
if $R$ is $\mathcal{S}$-Cohen-Macaulay over itself.
\end{definition}

In the sequel, we need the following theorem, which provides some
equivalent conditions to Definition 3.4.

\begin{theorem}\label{pro}
Let $(R,\fm)$ be a local ring and $M$ a finitely generated
$R$-module. Then the following are equivalent:
\begin{enumerate}\item[(i)]$M$ is  an $\mathcal{S}$-Cohen-Macaulay
$R$-module.
\item[(ii)]For any $\fp\in\mathcal{S}-\Supp_{R}(M)$, $M_{\fp}$ is
Cohen-Macaulay and $\Ht_M(\fp)+\dim( R/\fp)=\dim M$.

\item[(iii)] For any $\fp\in\mathcal{S}-\Supp_{R}(M)$,
$E.\grade_{\mathcal{S}}(\fp,M)
=\mathcal{S}-\Ht_{M}(\fp)=\Ht_{M}(\fp)$ and $\Ht_M(\fp)+\dim
(R/\fp)=\dim M$.

\item[(iv)]$\depth_{R_{\fp}}(M_{\fp})=\dim M-\dim( R/\fp)$, for any
$\fp\in\mathcal{S}-\Supp_{R}(M)$.

\item[(v)]If $x_1,\cdots,x_{d}$ is a system of parameters
for $M$, then $\dim (R/\fp)=\dim M-i$ for all
$\fp\in\mathcal{S}-\Ass_{R}(M/(x_1,\cdots,x_{i})M)$ and all $0\leq
i\leq d:=\dim M$.

\item[(vi)]
$\Ht(\fp/\fq)+\Ht_{M}(\fq)=\Ht_{M}(\fp)$ for all $\fq\subseteq\fp$
of $\mathcal{S}-\Supp_{R}(M)\cup\{\fm\}$, $M_{\fp}$ is Cohen
Macaulay for all $\fp\in\mathcal{S}-\Supp_{R}(M)$ and $\dim(
R/\fp)=\dim M$ for all $\fp\in\min(\mathcal{S}-\Supp_{R}(M))$.
\end{enumerate}
\end{theorem}

{\bf Proof.} $(ii)\Leftrightarrow (iv)$ is easy and
$(ii)\Leftrightarrow (iii)$ follows from Proposition
\ref{localization}.

$(i)\Rightarrow (iii)$ Let $\fp\in\mathcal{S}-\Supp_{R}(M)$ and
suppose that $\dim (R/\fp)=\dim M-i$. Then there exists a subset of
system of parameters $x_1,\cdots,x_{i}$ for $M$ that belongs to
$\fp$. In view of our assumption, the sequence $x_1,\cdots,x_{i}$ is
a weak $M$-sequence with respect to $\mathcal{S}$. So by Lemma 3.2
and Theorem \ref{comm}(ii), $$\mathcal{S}-E.\grade_{R}(\fp,M)\leq
\Ht_{M}(\fp)\leq i\leq \mathcal{S}-E.\grade_{R}(\fp,M).$$Therefore
$\mathcal{S}-E.\grade_{R}(\fp,M)= \Ht_{M}(\fp)=i$.

$(iii)\Rightarrow (v)$  First we claim that if
$\underline{x}:=x_1,\cdots,x_{\ell}$ is a subset of a system of
parameters for $M$, then $\underline{x}$ form a weak $M$-sequence
with respect to $\mathcal{S}$. We prove this claim by induction on
$\ell$. For the case $\ell=0$, we have nothing to prove. Now suppose
inductively, $\ell>0$ and the result has been proved for all subset
of system of parameters of length less than $\ell$.  Then, by the
inductive hypothesis $x_1,\cdots,x_{\ell-1}$ is a weak $M$-sequence
with respect to $\mathcal{S}$. Let
$\fp\in\mathcal{S}-\Ass_{R}(M/(x_1,\cdots,x_{\ell -1})M)$. In view
of Lemma 2.3, we have that $x_1,\cdots,x_{\ell-1}$ is a maximal weak
$M$-sequence with respect to $\mathcal{S}$ in $\fp$. By the
assumption and Theorem \ref{comm}(i), we have $\dim (R/\fp)=\dim
M-\ell+1$. Since $\underline{x}$ is a subset of a system of
parameters for $M$, one can deduce that $x_{\ell}\notin\fp$. Now the
implication $(i)\Rightarrow (ii)$ of Lemma 2.3 shows that
$\underline{x}$ is a weak $M$-sequence with respect to
$\mathcal{S}$, as desired claim.

Let $\underline{x}:=x_1,\cdots,x_{i}$ be a subset of a system of
parameters for $M$ and let
$\fp\in\mathcal{S}-\Ass_{R}(\frac{M}{(x_1,\cdots,x_{i})M})$. So
$x_1,\cdots,x_{i}$ is a maximal weak $M$-sequence with respect to
$\mathcal{S}$ in $\fp$ and consequently
$i=\mathcal{S}-E.\grade_{R}(\fp,M)$. Now the assertion becomes
clear.

$(v)\Rightarrow (i)$ Let $\underline{x}:=x_1,\cdots,x_{d}$ be a
system of parameters for $M$. In view of our assumption, we have
that $x_{i}\notin \bigcup_{\fp\in
\mathcal{S}-\Ass_{R}M/(x_{1},\cdots,x_{i-1})M}\fp$ for all $i=1,
\cdots, d$. Therefore this implication follows from Lemma 2.3 and
Definition 3.4.

$(ii)\Rightarrow (vi)$ Assume that
$\fp\in\min(\mathcal{S}-\Supp_{R}(M))$. Let $\fq\in\Supp M $ be such
that $\fq \subseteq \fp$. The epimorphism $R/\fq\longrightarrow
R/\fp$ shows that $R/\fq\notin\mathcal{S}$ and consequently
$\fp=\fq$.  Hence $\Ht_M(\fp)= 0$.  Therefore, by assumption we get
that $\dim R/\fp = \dim M$.

Let $\fp,\fq\in\mathcal{S}-\Supp_{R}(M)\cup\{\fm\}$ be such that
$\fq\subseteq\fp$. In order to show that
$\Ht_M(\fp)=\Ht(\fp/\fq)+\Ht_{M}(\fq)$, we can assume that
$\fp\neq\fq$. To do this, first we assume that $\fp\neq\fm$. So
$\fp\in\mathcal{S}-\Supp_{R}(M)$. Since $M_{\fp}$ is Cohen-Macaulay,
we have:
$$\Ht_M(\fp)= \dim (M_{\fp})= \dim (R_{\fp}/\fq R_{\fp})+ \Ht_{M_{\fp}}(\fq
R_{\fp})=\Ht(\fp/\fq)+\Ht_{M}(\fq).$$ Now, we consider the case
$\fp=\fm$. Note that $\fp\neq\fq$. Therefore
$\fq\in\mathcal{S}-\Supp_{R}(M)$. So
$$
\Ht_{M}(\fq) + \Ht(\fm/\fq)=\dim (M_{\fq})+\dim (R/\fq)=\dim
M=\Ht_{M}(\fm).
$$
It remains to show the implication $(vi)\Rightarrow (ii)$. This is
clear. $\Box$

For an $R$-module $L$, we denote $\{\fp \in
\mathcal{S}-\Supp_{R}L|\dim L=\dim (R/ \fp)\}$ by
$\mathcal{S}-\Assh_R L$. Now, we are ready to present some of the
basic properties of $\mathcal{S}$-Cohen-Macaulay modules.

\begin{proposition}
Let $(R,\fm)$ be a local ring, $M$ a finitely generated
$\mathcal{S}$-Cohen-Macaulay and $x$  a weak $M$-sequence in $\fm$
with respect to $\mathcal{S}$. If $\mathcal{S}-\Assh_R (M/xM)\neq
\emptyset$, then $M/xM$ is an $\mathcal{S}$-Cohen-Macaulay
$R$-module.
\end{proposition}

{\bf Proof.} First, we show that $\dim (M/xM)=\dim M -1$. It is
known that $\dim M\leq\dim (M/xM)+1$. Let $\fq\in
\mathcal{S}-\Assh_{R}(M/xM)$. Then by Theorem 2.8(i), Lemma 3.2 and
Theorem 3.5(ii) we have:$$
\begin{array}{ll}
1&\leq \mathcal{S}-E.\grade _{R}(\fq,M)\\&\leq
\mathcal{S}-\Ht_{M}(\fq)\\&=\Ht_{M}(\fq)\\&=\dim M-\dim(R/
\fq)\\&=\dim M-\dim (M/xM)\leq 1,
\end{array}
$$
which implies that $\dim (M/xM)=\dim M -1$.

Let $\fp\in\mathcal{S}-\Supp_{R}(M/xM)$. So $x\in\fp$ and
$\fp\in\mathcal{S}-\Supp_{R}(M)$. By Theorem \ref{pro}(ii) $M_{\fp}$
is Cohen-Macaulay and $\Ht_M(\fp)+\dim(R/\fp)=\dim M$. Also by Lemma
2.3 the element $x/1$ of the local ring $R_{\fp}$ becomes  an
$M_{\fp}$-sequence. This implies that $M_{\fp}/xM_{\fp}$ is
Cohen-Macaulay and $\Ht_{\frac{M}{xM}}(\fp)=\Ht_{M}(\fp)-1$.
Therefore $\Ht_{\frac{M}{xM}}(\fp)+\dim( R/\fp)=\dim (M/xM)$. Again
by Theorem \ref{pro}(ii), $M/xM$ is $\mathcal{S}$-Cohen-Macaulay.
$\Box$

\begin{corollary}
Let $(R,\fm)$ be a local ring, $x$ a weak $M$-sequence in $\fm$ with
respect to $\mathcal{S}$ and $M$ a finitely generated $R$-module.
Then the following hold:
\begin{enumerate}
\item[(i)] If $\frac{M}{xM}$ is  equidimensional and $M$ is
$\mathcal{S}$-Cohen-Macaulay, then $\frac{M}{xM}$ is
$\mathcal{S}$-Cohen-Macaulay.

\item[(ii)] If $M$ is an f-module (generalized f-module), then $\frac{M}{xM}$ is
an
f-module(generalized f-module).
\end{enumerate}
\end{corollary}

{\bf Proof.} $(i)$ Without loss of generality, we can assume that
$M/xM\notin\mathcal{S}$. It follows  from Lemma \ref{key} that there
exists $\fq\in\mathcal{S}-\Supp_R (M/xM)$. Let $\fp\subseteq\fq$
such that $\fp\in\min(\Supp_R (M/xM))$. The natural epimorphism
$R/\fp\longrightarrow R/\fq$ shows that $R/\fp\notin\mathcal{S}$. On
the other hand  $\frac{M}{xM}$ is an equidimensional $R$-module.
Therefore $\mathcal{S}-\Assh_R (M/xM)\neq\emptyset$.

$(ii)$ Set $\mathcal{S}:=\{N\in R\textbf{-Mod}|\Supp_R N \subseteq
\{\fm \}\}$ ($\mathcal{S}:=\{N\in R\textbf{-Mod}|\dim N\leq 1\}$).
Then the notions of $\mathcal{S}$-Cohen-Macaulay and  f-module
(generalized f-module) are equivalent. Now, we prove the claim. We
can assume that $M/xM\notin\mathcal{S}$. Let $\fp\in \Assh_R
(M/xM)$. Therefore $\dim R/\fp>0$ ($\dim R/\fp>1$). Thus in both
cases we have $\mathcal{S}-\Assh_R (M/xM)\neq\emptyset$. $\Box$

\begin{remark} Let $(R,\fm)$ be a local ring and $M$ a finitely generated
$R$-module. Assume that $x$ is an $M$-regular element. This is well
known that $M$ is Cohen-Macaulay if and only if $M/xM$ is
Cohen-Macaulay. Having this fact and Corollary 3.7 in mind, one
might ask whether the converse of Corollary 3.7(i) is true. This is
not necessary true. To see this, consider the following example. Let
$(R,\fm)$ be a 2-dimensional local domain such that its completion
$\widehat{R}$ has an associated prime $\fp$ with $\dim
\widehat{R}/\fp=1$. Such rings constructed by Nagata \cite[page 203,
Example 2]{Na}. Theorem \ref{pro}(iv) implies that $\widehat{R}$ is
not an f-module. Since $\depth \widehat{R}=\depth R>0$, there exists
an element  $x\in\widehat{R}$, which is $\widehat{R}$-regular. So
$\dim \widehat{R}/x\widehat{R}=1$. On the other hand, any
1-dimensional local ring is an f-module. Therefore
$\widehat{R}/x\widehat{R}$ is an f-module.
\end{remark}

Let $f:R\longrightarrow A$ be a flat homomorphism of rings and
$\mathcal{S}$ a Serre class of $A$-modules. Set
$\mathcal{S}^c=\{M\in R\textbf{-Mod}|M\otimes_R A\in\mathcal{S}\}$.
It is routine to show that $\mathcal{S}^c$ is a Serre class  of
$R$-modules.

\begin{theorem} Let $f:(R,\fm)\longrightarrow (A,\fn)$ be a
flat local  homomorphism of Noeherian local rings. Let $\mathcal{S}$
be a Serre class of $A$-modules and $M$  a finitely generated
$R$-module. Then $M\otimes_R A$ is $\mathcal{S}$-Cohen-Macaulay if
and only if the following three conditions are satisfied:
\begin{enumerate}
\item[(i)]
$M$ is $\mathcal{S}^c$-Cohen-Macaulay.
\item[(ii)] $\frac{A_{\fq}}{f^{-1}(\fq)A_{\fq}}$ is
Cohen-Macaulay, for all $\fq\in\mathcal{S}-\Supp_{A}(M\otimes_R A)$.
\item[(iii)]$\Ht(\frac{\fq}{f^{-1}(\fq) A})+\dim (\frac{A}{\fq})=\dim
(\frac{A}{f^{-1}(\fq) A})$,  for all
$\fq\in\mathcal{S}-\Supp_{A}(M\otimes_R A)$.
\end{enumerate}
\end{theorem}

{\bf Proof.} We first bring the following easy results (a) and (b).

(a) Let $\fp\in\mathcal{S}^c-\Supp_{R}(M)$. The condition $R/
\fp\notin\mathcal{S}^c$ implies that $A/  \fp A\notin\mathcal{S}$.
In view of Lemma \ref{key}, one can find
$\fq'\in\mathcal{S}-\Supp_{A}(A/\fp A)$. Let $\fq$ be a minimal
prime ideal of $\fp A$ such that $\fq\subseteq\fq'$. Hence
$\fp=f^{-1}(\fq)$.  Also the natural epimorphism
$A/\fq\longrightarrow A/\fq'$ shows that $A/\fq\notin\mathcal{S}$.
Therefore $\fq\in\mathcal{S}-\Supp_{A}(M\otimes_R A)$.

(b) Let the situation and notation  be as in (a). Then we have
$$
\begin{array}{ll}
\Ht_M(\fp)+\dim( R/\fp)&=\dim(A_{\fq}/\fp A_{\fq})+\dim
M_{\fp}+\dim( R/\fp)\\&= \dim (M_{\fp}\otimes_{R_{\fp}}A_{\fq})+\dim
(R/ \fp)\\&= \Ht_{M\otimes_R A}(\fq)+\dim (R/ \fp)\\&=
\Ht_{M\otimes_R A}(\fq)+\dim (A/\fp A )-\dim (A/ \fm A)\\&\geq
\Ht_{M\otimes_R A}(\fq)+\dim (A/\fq ) -\dim (A/ \fm A),
\end{array}
$$
where the second equality follows from \cite[Theorem A.11(ii)]{BH},
for the natural flat homomorphism $R_{\fp}\longrightarrow A_{\fq}$,
and the forth equality follows from  \cite[Theorem A.11(i)]{BH}, for
the natural flat homomorphism $R/ \fp \longrightarrow A/\fp A$.

Now, we are ready to prove our claims.

Assume that $M\otimes_R A$ is $\mathcal{S}$-Cohen-Macaulay. Let
$\fp\in\mathcal{S}^c-\Supp_{R}(M)$. Then by (a) there exists $\fq\in
\mathcal{S}-\Supp_{A}(M\otimes_R A)$ such that $\fp=f^{-1}(\fq)$ and
$\fq$ is a minimal prime ideal of $\fp A$. By our assumption and
Theorem \ref{pro} $(M\otimes_R A)_{\fq}\cong
M_{\fp}\otimes_{R_{\fp}} A_{\fq}$ is Cohen-Macaulay. So $M_{\fp}$
and $A_{\fq}/\fp A_{\fq}$ are Cohen Macaulay.

On the other hand from (b) and Theorem \ref{pro}(ii) we have $$
\begin{array}{ll}
\Ht_M(\fp)+\dim (R/\fp)&\geq \Ht_{M\otimes_R A}(\fq)+\dim (A/\fq)
-\dim (A/ \fm A)\\&=\dim (M\otimes_R A)-\dim (A/ \fm A)\\&=\dim M.
\end{array}
$$
Consequently $\Ht_M(\fp)+\dim (R/\fp)=\dim M$. Again by Theorem
\ref{pro}(ii),  $M$ becomes  $\mathcal{S}^c$-Cohen-Macaulay. This
completes the proof of (i).

With the notation as (ii), the desired result of (ii) follows from
the Cohen-Macaulayness of $(M\otimes_R A)_{\fq}\cong
M_{f^{-1}(\fq)}\otimes_{R_{\fp}} A_{\fq}$.

In order to prove (iii) first we claim that:

(c) Let $\fp\in\mathcal{S}^c-\Supp_{R}(M)$. If $Q$ is a minimal
prime ideal of $\fp A$ and $A/Q\notin\mathcal{S}$, then $\dim (A/Q)=
\dim(A/ \fp A)$.

To establish the claim, consider the following:
$$
\begin{array}{ll}
\dim(A/ \fp A)&\geq \dim(A/ Q)\\&= \Ht_{M\otimes_R A}(\fn)-\Ht
_{M\otimes_R A}(Q)\\&= \Ht_{M\otimes_R A}(\fn)-\Ht _{M}(\fp)\\&=\dim
M+ \dim (A/\fm A)-\Ht _{M}(\fp)\\&=\dim (R/ \fp)+\dim (A/ \fm
A)\\&=\dim(A/ \fp A), \end{array}
$$
where the first equality follows from Theorem \ref{pro}(ii) and the
second equality follows from \cite[Theorem A.11(ii)]{BH}, for the
natural flat homomorphism $R_{\fp}\longrightarrow A_{Q}$. Note that
we have $Q\in\mathcal{S}-\Supp_{A}(M\otimes_R A)$ and $Q\cap R=\fp$.

Let $\fq\in\mathcal{S}-\Supp_{A}(M\otimes_R A)$ and
$\fp=f^{-1}(\fq)$. Then $A/\fp A\notin \mathcal{S}$. So by (a) we
can find $\fq_1\in\mathcal{S}-\Supp_{A}(M\otimes_R A)$  such that
$\fq_1$ is a minimal prime ideal of $\fp A$ and $\fq_1\subseteq
\fq$. Keep in mind that $M\otimes_R A$ is
$\mathcal{S}$-Cohen-Macaulay. Hence Theorem \ref{pro}(vi) implies
that
$$
\begin{array}{ll}
\dim (A/\fq)+\Ht(\fq/\fp A)&\geq\Ht(\fn/\fq)+\Ht(\fq/\fq _1)\\&=
\Ht_{M\otimes_R A}(\fn)-\Ht_{M\otimes_R A}(\fq)+\Ht_{M\otimes_R
A}(\fq)-\Ht_{M\otimes_R A}(\fq_1)\\&=\Ht(\fn/\fq_1)
\\&=\dim (A/\fq_1)\\&=\dim (A/\fp A),
\end{array}
$$
where the last equality follows from (c). Consequently $\Ht(\fq/\fp
A)+\dim (A/\fq)=\dim (A/\fp A)$. This completes the proof of (iii).

Conversely,  assume that the conditions (i), (ii) and (iii) hold. In
order to show that $M\otimes_R A$ is $\mathcal{S}$-Cohen-Macaulay,
in view of Theorem \ref{pro}(ii), we need to prove the following
claims:

\begin{enumerate}
\item[(d)] $(M\otimes_R A)_{\fq}$ is Cohen-Macaulay, for all
$\fq\in\mathcal{S}-\Supp_{A}(M\otimes_R A)$.
\item[(e)]$\Ht_{M\otimes_R A}(\fq)+\dim (A/\fq)=\dim (M\otimes_R A)$, for all
$\fq\in\mathcal{S}-\Supp_{A}(M\otimes_R A)$.
\end{enumerate}

The claim (d) follows immediately from (i) and (ii). Now, we shall
achieve (e). Let $\fq\in\mathcal{S}-\Supp_{A}(M\otimes_R A)$. Set
$\fp:=f^{-1}(\fq)$. Then we have
$$
\begin{array}{ll}
\dim (A/\fq)+\Ht_{M\otimes_R A}(\fq)&=\dim (A/\fq )
+\Ht_M(\fp)+\Ht(\fq/\fp A )\\&=\dim (A/\fq ) +\dim M-\dim (R/
\fp)+\Ht(\fq/\fp A)\\&=\dim (A/\fq ) +\dim (M\otimes_R A)-\dim (A/
\fm A)-\dim (R/ \fp)+\Ht(\fq/\fp A)
\\&=\dim (A/\fq )+\dim (M\otimes_R A)-\dim (A/ \fp A)+\Ht(\fq/\fp A)
\\&=\dim (M\otimes_R A),
\end{array}
$$
where the last equality follows from (iii). $\Box$

For the ring extension $(R,\fm)\longrightarrow (\widehat{R},\fm
{\widehat{R}})$, Theorem 3.9 becomes  much simpler:

\begin{proposition}
Let $(R,\fm)$ be a local ring and $M$  a finitely generated
$R$-module.  If $\mathcal{S}$ is a  Serre class of
$\widehat{R}$-modules and $M\otimes_R \widehat{R}$ is
$\mathcal{S}$-Cohen-Macaulay, then $M$ is $\mathcal{S}^c$-Cohen
Macaulay.
\end{proposition}

{\bf Proof.}  Let $\underline{x}:=x_1,\cdots,x_{\ell}$ be a system
of parameters for $M$.  So $\underline{x}$ is a system of parameters
for $M\otimes_R \widehat{R}$. Therefore, $((x_{1},\dots,
x_{i-1})M\otimes_R \widehat{R}:_{M\otimes_R \widehat{R}}
x_{i})/(x_{1},\dots, x_{i-1})(M\otimes_R
\widehat{R})\in\mathcal{S}$. The isomorphism
$$\frac{((x_{1},\dots, x_{i-1})M\otimes_R \widehat{R}:_{M\otimes_R \widehat{R}}
x_{i})}{(x_{1},\dots, x_{i-1})(M\otimes_R \widehat{R})}\cong
\frac{((x_{1},\dots, x_{i-1})M:_{M} x_{i})}{(x_{1},\dots,
x_{i-1})M}\otimes_R \widehat{R}$$ completes the proof.
$\Box$\\

Let $(R,\fm)$ be a local ring and $M$ a $d$-dimensional finitely
generated $R$-module.  If $d=0$ or $M=0$, we set $\fa(M)=R$. In the
other case, we set $\fa(M):=\fa_0(M). \cdots. \fa_{d-1}(M)$, where
$\fa_i(M)= \Ann_R(
 H^i_{\fm}(M))$ for all $i=0,\cdots,d-1$.

\begin{theorem}
Let $(R,\fm)$ be a local ring and $M$ a $d$-dimensional finitely
generated $R$-module. Consider the following conditions:
\begin{enumerate}\item[(i)] $R/\fa(M)\in\mathcal{S}$.
\item[(ii)] $M$ is  an $\mathcal{S}$-Cohen-Macaulay $R$-module.
\end{enumerate}
Always  $(i)$ implies $(ii)$ and if $R$ is a quotient of a
Gorenstein local ring, then $(ii)$ implies $(i)$.
\end{theorem}

{\bf Proof.} In the cases $\dim M=0$ and $M=0$, the desired claims
are trivial. Therefore without loss of generality we can assume that
$\dim M>0$.

$(i)\Rightarrow (ii)$ Let $\underline{x}:=x_1,\cdots,x_{d}$ be a
system of parameters for $M$. Set $I:=\underline{x}R$ and consider
$$\tau_{\underline{x}}(M):=\bigcap \Ann_R(\frac{((x_{1}^t,\dots,
x_{i-1}^t) M:_{M} x_{i}^t)}{(x_{1}^t,\dots, x_{i-1}^t)M}),$$ where
the intersection is taken over all $t\in\mathbb{N}$ and $1\leq i\leq
d$. We denote $\Ann_R (H^i_{I}(M))$ by $\fa^i_I(M)$. By
\cite[Theorem 3(a)]{Sch}, we have $\fa^0_I(M). \cdots .
\fa^{d-1}_I(M)\subseteq \tau_{\underline{x}}(M)$. Note that
$\sqrt{I+\Ann_R M}=\fm$, since $\underline{x}$ is a system of
parameters for $M$. The independence theorem for local cohomology
modules implies that $H^i_{\fm}(M)\cong H^i_{I}(M)$ and consequently
$\fa^i_I(M)= \fa_i(M)$. Therefore $\fa(M)\subseteq
\Ann_R(\frac{((x_{1},\dots, x_{i-1})M :_{M} x_{i})}{(x_{1},\dots,
x_{i-1})M})$ for all $1\leq i\leq d$. The condition
$R/\fa(M)\in\mathcal{S}$ implies that $R/\Ann_R(\frac{((x_{1},\dots,
x_{i-1})M :_{M} x_{i})}{(x_{1},\dots, x_{i-1})M})\in\mathcal{S}$ for
all $1\leq i\leq d$. On the other hand the $R$-modules
$R/\Ann_R(\frac{((x_{1},\dots, x_{i-1})M :_{M} x_{i})}{(x_{1},\dots,
x_{i-1})M})$ and $\frac{((x_{1},\dots, x_{i-1})M :_{M}
x_{i})}{(x_{1},\dots, x_{i-1})M}$ have same supports for all $1\leq
i\leq d$. Therefore, it is enough to apply Lemma \ref{key} to obtain
that $\underline{x}$ is a weak $M$-sequence with respect to
$\mathcal{S}$.\\

$(ii)\Rightarrow (i)$  Suppose the contrary that,
$R/\fa(M)\notin\mathcal{S}$. In view of Lemma \ref{key} we can find
$\fp\in\V(\fa(M))$ such that $R/\fp\notin\mathcal{S}$. Next we show
that $\depth_{R_{\fp}} (M_{\fp})+\dim (R / \fp) <\dim(M)$. So in
view of Theorem \ref{pro}(ii), we get a contradiction.

Now, we do this. We have $\fp\in\V(\fa_i(M))$, for some $i<\dim M$.
Let $(R',\fn)$ be a Gorenstein ring of dimension $r'$ for which
there exists a surjective ring homomorphism $f:R' \longrightarrow
R$. The local duality theorem \cite[Theorem 11.2.6]{BS} implies that
$$H^{i}_{\fm}(M)\cong
\Hom_{R}(\Ext^{r'-i}_{R'}(M,R'),E_{R}(R/\fm)),$$ where
$E_{R}(R/\fm)$ is the injective envelope of $R/\fm$. Hence
$\fp\in\Supp_R(\Ext^{r'-i}_{R'}(M,R'))$. Let $t:=\dim R/\fp$ and
$\fp':=f^{-1}(\fp)$. Now $R'_{\fp'}$ is a Gorenstein local ring and
$t=\dim (R'/\fp')$. Since $R'$ is a Gorenstein ring, so $\dim
R'_{\fp'}=\dim R'-\dim (R'/\fp')=r'-t$.

Let $f':R'_{\fp'} \longrightarrow R_{\fp}$ be  the surjective ring
homomorphism, which induced by $f$. There is an
$R_{\fp}$-isomorphism,
$\Ext^{r-i}_{R'_{\fp'}}(M_{\fp},R'_{\fp'})\cong
(\Ext^{r-i}_{R'}(M,R'))_{\fp}$. Again the local duality theorem
implies that $$H^{i-t}_{\fp R_{\fp}}(M_{\fp})\cong
\Hom_{R_{\fp}}(\Ext^{r'-i}_{R'_{\fp'}}(M_{\fp},R'_{\fp'}),E_{R_{\fp}}
(R_{\fp}/\fp R_{\fp})),$$ as $R_{\fp}$-module. Therefore
$H^{i-t}_{\fp R_{\fp}}(M_{\fp})\neq0$. Consequently
$$\depth_{R_{\fp}}(M_{\fp})\leq i-t<\dim M-t=\dim M-\dim (R / \fp). \ \ \Box$$

The following example shows that the assumption, $R$ is a quotient
of a Gorenstein local ring,  in  Theorem 3.11  is needed.

\begin{example} By \cite[Example 3.4]{NM}, there exists a 3-dimensional local ring
$(R,\fm)$ such that $\dim (R/\fa(R))=3$. Denote  the class of all
finitely generated $R$-modules of Krull dimension less than 2, by
$\mathcal{S}$. It is a Serre class of $R$-modules. The assumption
$\dim (R/\fa(R))=3$ implies that $R/\fa(R)\notin\mathcal{S}$. Also,
the example shows that $R$ is a generalized f-ring, i.e.
$\mathcal{S}$-Cohen-Macaulay.
\end{example}

Denote by $\mathcal{NCM}(M)$ the non Cohen-Macaulay locus of M, i.e.
$$ \mathcal{NCM}(M)=\{\fp\in \Spec R | M_{\fp} \textit{ is not Cohen-Macaulay}\}.$$
Assume that  $R$ is a quotient of a Cohen-Macaulay ring. It is well
known that the non Cohen-Macaulay locus of M is a closed subset of
$\Spec R$ with respect to the Zariski topology i.e.
$\mathcal{NCM}(M)=V(\fa_M)$, for some ideal $\fa_M$ of $R$.
Therefore such ideals are unique up to radical.

\begin{theorem}
Let $(R,\fm)$ be a local ring which is a quotient of a
Cohen-Macaulay ring and $M$ a  finitely generated $R$-module. Then
the following are equivalent:
\begin{enumerate}\item[(i)] $M$ is  an $\mathcal{S}$-Cohen-Macaulay
$R$-module.
\item[(ii)]$R/\fa_M\in\mathcal{S}$ and $\dim (R/\fp)=\dim M$ for all prime ideals
$\fp\in\min (\mathcal{S}-\Supp_{R}(M))$.
\end{enumerate}
\end{theorem}

{\bf Proof.} $(i)\Rightarrow (ii)$ In order to prove
$R/\fa_M\in\mathcal{S}$, it is enough to show that
$R/\fp\in\mathcal{S}$ for all $\fp\in\V(\fa_M)$. Let
$\fp\in\V(\fa_M)$. Since $M_{\fp}$ is not Cohen-Macaulay, Theorem
\ref{pro}(ii) implies that $R/\fp\in\mathcal{S}$.

Let $\fp\in\min(\mathcal{S}-\Supp_{R}(M))$. Then by the implication
$(i)\Rightarrow (vi)$ of Theorem  \ref{pro}, we have $\dim
(R/\fp)=\dim M$.

$(ii)\Rightarrow (i)$  Let $\fp\in\mathcal{S}-\Supp_{R}(M)$. In
particular $R/ \fp \notin\mathcal{S}$. In view of Lemma \ref{key}
and $R/\fa_M\in\mathcal{S}$, we have $\fp\notin\V(\fa_M)$. So
$M_{\fp}$ is a Cohen-Macaulay $R_{\fp}$-module.

Set $t:=\dim (R/\fp)$ and $s:=\Ht_M (\fp)$. Hence we have following
saturated chains of prime
ideals$$\fp=\fp_0\subset\cdots\subset\fp_t=\fm,$$
$$\fq=\fq_0\subset\cdots\subset\fq_s=\fp.$$
By concatenating these chains, we get the following saturated chain
of prime ideals:
$$\fq=\fq_0\subset\cdots\subset\fq_s=\fp=\fp_0\subset\cdots\subset\fp_t=\fm.$$
The epimorphism $R/\fq_0\longrightarrow R/\fp$ shows that
$R/\fq_0\notin\mathcal{S}$ and consequently $\fq_0\in\min(
\mathcal{S}-\Supp_{R}(M))$. So $\dim (R/\fq_0)=\dim M$. On the other
hand, $R$ is a quotient of a Cohen-Macaulay ring. This implies that
$\Supp (R/\fq_0)$ is  catenary. Therefore $$\dim M=\dim
(R/\fq_0)=s+t=\Ht_M (\fp)+\dim (R/\fp).$$ Consequently, Theorem
\ref{pro}(ii) implies that $M$ is an $\mathcal{S}$-Cohen-Macaulay
$R$-module. $\Box$



\begin{thebibliography}{99}

\bibitem[A]{A}{K. Ahmadi-Amoli}, {\it
Local cohomology and cohomology of certain generalized Hughes
complexes}, Comm. Algebra , {\bf 26}, (4), (1998),  1305--1317.

\bibitem[BH]{BH}{W. Bruns and J. Herzog}, {\it Cohen-Macaulay rings},
Cambridge Univ. Press, {\bf 39}, Cambridge, (1998).

\bibitem[BN]{BN}{ M. Brodmann}, {L. T. Nhan}, {\textit A finiteness
result for associated primes of certain Ext-modules}, To appear in
Comm. Algebra.

\bibitem[BS]{BS}{ M. Brodmann}, { R. Y. Sharp}, {\it
Local cohomology, An algebraic introduction with geometric
applications}, Cambridge Univ. Press, {\bf 60}, (1998).

\bibitem[CHK]{CHK}{N. T. Cuong}, {N. V. Hoang}, {P. H. Khanh},
{\textit Asymptotic stability of certain sets of assuciated prime
ideals of local cohomology modules}, arXiv:0804.0964v1.

\bibitem[CST]{CST}{ N. T. Cuong}, { P. Schenzel}, {N. V. Trung},  {\it
Verallgemeinerte Cohen–Macaulay moduln}, Math. Nachr,  {\bf 85},
(1978), 57–-73.

\bibitem[GR]{GR}{O . Gabber}, { L. Ramero}, {\it Almost ring
theory}, Springer LNM, {\bf1800},  (2003).



\bibitem[LT]{LT}{R. Lu}, {Z. Tang}, {\it The f-depth
of an ideal on a module}, Proc. Amer.
Math. Soc., {\bf 130}, (2002), 1905--1912.




\bibitem[Ma]{Ma}{ A. Mafi},  {\it Some results on local cohomology modules},
Arch. Math., {\bf 87}, (3), (2006), 211--216.

\bibitem[Mel1]{Mel1}{ L. Melkersson},  {\it Modules cofinite
with respect to an ideal}, J. Algebra, {\bf 285}, (2005), 649–-668.


\bibitem[Mel2]{Mel2}{ L. Melkersson},  {\it Some applications
of a criterion for Artinianness of a module}, J. Pure Appl. Algebra
{\bf 101}, (1995), 291-303.


\bibitem[Na]{Na}{M. Nagata}, {\it Local rings}. Interscience, (1969).

\bibitem[N]{N}{L. T. Nhan}, {\it On generalized regular sequences and the
finiteness for associated primes of local cohomology modules}. Comm.
Algebra {\bf 33}, (3), (2005),  793-806.

\bibitem[NM]{NM}{L. T. Nhan}, {M. Morales}, {\it Generalized
f-modules and the associated primes of local cohomology modules}.
Comm. Algebra {\bf 34},  (2006),  863--878.

\bibitem[Q]{Q}{ N. Quoc Chinh},  {\it
On reducing sequences and an application to local cohomology
modules},  Bull. Iranian Math. Soc., {\bf 31}, (2), (2005), 1--11.

\bibitem[RSS]{RSS}{ P. Roberts}, { A. K. Singh}, { V. Srinivas}, {\it Annihilators of local
cohomology in characteristic zero}, Illinois J. of Math., {\bf51},
(1), (2007),  237-–254.


\bibitem[Sch]{Sch}{P. Schenzel},  {\it
Cohomological annihilators}, Proc. Phil. Soc., {\bf 91}, (1982),
345--350.


\bibitem[St]{St}{B. Stenstr\"{o}m},  {\it
Rings of quotients}, Die Grundlehren der Mathematischen
Wissenschaften {\bf217}, Springer-Verlag, New York-Heidelberg,
(1975).

\bibitem[Str]{Str}{ Strooker},  {\it
Homological questions in local algebra}, London mathematical Lecture
Note series {\bf 145}, (1990).

\bibitem[SV]{SV}{J. St$\ddot{u}$ckrad}, {W. Vogel},
{\it Buchsbaum rings and applications}, Berlin: Web Deutsecher
Verlag der Wissenschaften, (1986).


\bibitem[T]{T}{R. Takahashi},  {\it Classifying subcategories of modules over a commutative noetherian ring},  preprint
(2006).

\end{thebibliography}
\end{document}